\documentclass[12pt]{article}%
\usepackage{amsmath}
\usepackage{graphicx}
\usepackage{amsfonts}
\usepackage{amssymb}
\usepackage{graphicx}%
\providecommand{\U}[1]{\protect\rule{.1in}{.1in}}
\providecommand{\U}[1]{\protect\rule{.1in}{.1in}}
\usepackage{cite}

\begin{document}

\begin{center}
{\Large \textbf{The effect of finite duration sources on modes and  \\[10pt] generalization of the d'Alembert solution}}

\bigskip
\bigskip

\textsc{J.\ S.\ Ben-Benjamin${}^{1,a}$ and L.\ Cohen${}^{2}$}

\bigskip

\textit{${}^{1}$Great Wing Inc., 2090 Pacific Blvd, Atlantic Beach, NY 11509}

\textit{${}^{2}$Institute for Quantum Science and Engineering, Texas A\&M, College Station, TX 77840}

${}^{a}$corresponding author: j.s.ben-benjamin@greatwing.com

\bigskip

(\today)
\end{center}


\hskip                    .40cm

\vskip.14in \hskip.5cm\begin{minipage}[t]{15.05cm}
	\noindent{\small \textbf{Abstract}: 
We investigate the evolution of dispersive waves governed by linear wave
equations, where a finite duration source is applied.
The resulting wave may be viewed as the superposition of modes before the source is
turned on and after it is turned off.
We consider the problem of relating the modes after the source term is turned
off to the modes before the source term was turned on.
We obtain explicit formulas in both the wavenumber and position representations.
A number of special cases are considered.
Using the methods presented,
we obtain a generalization of the d'Alembert solution which applies to linear wave equations with constant coefficients.
}\end{minipage}\vskip.2truein
\bigskip

\bigskip


\vskip.3in


\section{Introduction}

For propagating dispersive waves governed by linear wave equations with
constant coefficients and without any source terms, the wave may be decomposed
into modes. However, if there is a source term of finite duration, then for
the time that the source term is on, the wave cannot be described by modes. In
particular suppose the governing equation is {\cite{morse2,pinch}}
\begin{equation}
\mathcal{L}\,u(x,t)=%
\begin{cases}
0 & \text{for }t<-T,\\
f(x,t) & \text{for }-T<t<0,\\
0 & \text{for }t>0.
\end{cases}
\label{eq8}%
\end{equation}
where $u(x,t)$ is the wave, $f(x,t)$ is the source term, and the linear wave equation, $\mathcal{L}$, is
given by
\begin{equation}
\mathcal{L=}\sum_{n=0}^{N_{x}}b_{n}\frac{\partial^{n}}{\partial x^{n}}%
-\sum_{n=0}^{N_{t}}a_{n}\frac{\partial^{n}}{\partial t^{n}}%
\end{equation}
That is, the wave evolves without the source term up to time $t=-T$, at which
time, a source term $f(x,t)$ is applied up to a time $t=0$, and is then turned
off. The resulting wave may be described by modes only before the source term
is turned on and after it is turned off. We consider the problem of relating
the modes after the source term is turned off to the modes before the source
term was turned on. That is, we investigate the effect of the finite duration
source term on modes.

In the next section we describe the method of modes and also discuss how to
obtain the modes given the initial conditions for the wave.

\section{Method of Modes}

For waves governed by a wave equation
\begin{equation}
\mathcal{L}\,u(x,t)=0\label{eq01}%
\end{equation}
where%
\begin{equation}
\mathcal{L}\mathcal{=}\sum_{n=0}^{N_{x}}b_{n}\frac{\partial^{n}}{\partial
x^{n}}-\sum_{n=0}^{N_{t}}a_{n}\frac{\partial^{n}}{\partial t^{n}}%
\end{equation}
one substitutes $u=e^{ikx-i\omega t}$ into Eq.\ (\ref{eq01})
\begin{equation}
\mathcal{L}\,e^{ikx-i\omega t}=0
\end{equation}
to obtain the dispersion relation \cite{graf,jack,light,morse,tolt,whit}
\begin{equation}
\sum_{n=0}^{N_{x}}b_{n}(ik)^{n}-\sum_{n=0}^{N_{t}}a_{n}(-i\omega
)^{n}=0\label{eq03}%
\end{equation}
One may solve Eq.\ (\ref{eq03}) for $\omega$ as a function of $k$, or
for $k$ as a function of $\omega$; the choice depends on the physical
situation. For example, if we pluck a string at a particular time and let go,
then one would solve for $\omega$ as a function of $k$. This corresponds to
when $u(x,0)$ is the initial condition. Alternatively, if we are at a fixed
position and create a pulse as a function of time, $u(0,t)$, then one would
solve for $k$ as a function of $\omega.$ In this paper, we deal with the
former case where we solve for $\omega$ as a function of $k$.


Generally speaking, there is more than one solution to the dispersion
relation. Assuming that there are $M$ solutions, we write
\begin{equation}
\omega=\omega_{\ell}(k),\hspace{0.5in}\mbox{for }\ell=1,2,\cdots\!,M;
\label{eq6}%
\end{equation}
each $\omega_{\ell}(k)$ may be complex. When the source term is zero, the wave
is then given by the sum of modes
\begin{equation}
u(x,t)\,=\sum_{\ell=1}^{M}u_{\ell}(x,t)
\end{equation}
where $u_{\ell}(x,t)$ are the modes; each corresponds to a solution of the
dispersion relation, $\omega_{\ell}(k)$. The properties and evolution of the
modes and methods for obtaining them are described in
Sec.\ \ref{sec:spatial_formulation},
where we also discuss how to obtain the initial conditions of the modes from
the initial conditions for the wave.

Before the source is turned on, the wave is the superposition of modes
\begin{equation}
u(x,t)\,=\sum_{\ell=1}^{M}u_{\ell}^{B}(x,t),\qquad\text{for }t<-T
\label{eq:20230309_105453}%
\end{equation}
where $u_{\ell}^{B}(x,t)$ are the mode \emph{before} the source is turned on.

While the source is on, we write the solution as%
\begin{equation}
u(x,t)=\mu(x,t)\qquad\text{for }-T<t<0\,\, \label{eq:20230309_105727}%
\end{equation}
After the source is turned off, we can again express the solution in terms of
modes
\begin{equation}
u(x,t)\,=\sum_{\ell=1}^{M}u_{\ell}^{A}(x,t)\hspace{0.2in}\text{for }t>0
\label{eq:20230309_105545}%
\end{equation}
where $u_{\ell}^{A}$ are the modes \emph{after} the source is tuned off.

As mentioned in the Introduction, in this paper, we address the question of
how modes are affected by the source term. That is, we consider the
relation of the modes after the source $u_{\ell}^{A}$ to the modes before the
source $u_{\ell}^{B}.$

\bigskip

\subsection{\textbf{Notation: }}

\noindent
All integrals range from $-\infty$ to $\infty$ unless otherwise noted.

\bigskip

\noindent
$u(x,t)\,$and $S(k,t)$ are the wave expressed in position space and wavenumber
space (which is discussed in the next section).

\bigskip

\noindent
$u_{\ell}^{B}$ and $u_{\ell}^{A}$ are the modes \emph{before} and \emph{after} the source, respectively.

\bigskip

\noindent
$S_{\ell}^{B}(x,t)$ and $S_{\ell}^{A}(x,t)$ are the modes \emph{before} and \emph{after} the
source expressed in wavenumber space, respectively.

\bigskip

\noindent
$f(x,t)$ and $F(k,t)$ are the source term, expressed in position and wavenumber space, respectively.

\bigskip

\noindent
$u_{p}(k,t)$ and $S_{p}(k,t)$ are particular solutions in the time interval
$-T<t<0$.
\begin{align}
\varphi(k,t)  &  =S(k,t)\hspace{0.2in}\hspace{0.2in}\text{for }-T<t<0\\
\mu(x,t)  &  =u(x,t)\hspace{0.2in}\hspace{0.2in}\text{for }-T<t<0
\end{align}


\noindent
As usual, expressions such as $\frac{{\partial}{}}{{\partial}{t}}S_{2}%
^{A}(k,T)\,$ mean $\left.  \frac{{\partial}{}}{{\partial}{t}}S_{2}%
^{A}(k,t)\right\vert _{t=T}$.

\bigskip

\noindent
For expressions such as
\begin{align}
\frac{{\partial}{}}{{\partial}{t}}e^{i\omega_{2} t}S_{1}^{B}(k,t=0)
%
\end  {align}
we mean that the differentiation is done on $e^{i\omega_{2}t}S_{1}^{B}(k,t)$ and after the differentiation $t$ is set to zero only on the term indicated by $t=0$.
%

\bigskip

\noindent
We define wavenumber operator, $\mathcal{K}$ in the position representation
by
\begin{equation}
\mathcal{K}=\frac{1}{i}\frac{\partial{}}{\partial x} \label{eq14}%
\end{equation}

\section{Wavenumber formulation}

The wavenumber representation makes it easier to relate the modes after the
source term is tuned off to the modes before the source term is turned on. The
fundamental reason is that in the wavenumber formulation the equation of
motion for modes is relatively simple as we now show. The wave in wavenumber space, $S(k,t)$,
is
\begin{equation}
S(k,t)\,=\,\frac{1}{\sqrt{2\pi}}\int\,u(x,t)e^{-ikx}\,dx
\label{eq:20231029_125224}%
\end{equation}
and each mode is defined by
\begin{equation}
u_{\ell}(x,t)\,=\,{\frac{1}{\sqrt{2\pi}}}\int S_{\ell}(k,0)\,\,e^{ikx-i\omega
_{\ell}(k)t}\,dk \label{eq:9}%
\end{equation}
where $S_{\ell}(k,0)$ is the initial wavenumber mode, obtained
from the spatial mode $u_{\ell}(x,0)$ by way of
\begin{equation}
S_{\ell}(k,0)={\frac{1}{\sqrt{2\pi}}}\,\int u_{\ell}(x,0)\,\,e^{-ikx}\,dx
\label{eq12}%
\end{equation}

If one defines the time-dependent wavenumber mode as \cite{loug2}%
\begin{equation}
S_{\ell}(k,t)=e^{-i\omega_{\ell}(k)t}S_{\ell}(k,0) \label{eq:20230309_105924}%
\end{equation}
then $u_{\ell}(x,t)$ and $S_{\ell}(k,t)$ form Fourier transform pairs between
$x$ and $k$ for all time%
\begin{align}
u_{\ell}(x,t)\,  &  =\,\frac{1}{\sqrt{2\pi}}\int S_{\ell}(k,t)\,e^{ikx}\,dk\\
S_{\ell}(k,t)\,  &  =\,{\frac{1}{\sqrt{2\pi}}}\int u_{\ell}(x,t)\,e^{-ikx}\,dx
\end{align}
Each of the modes in wavenumber space evolve in a simple manner, namely
Eq.\ (\ref{eq:20230309_105924}), and satisfy the equation of motion \cite{cohen09} 
\begin{equation}
i\frac{{\partial}{}}{{\partial}{t}}S_{\ell}(k,t)=\omega_{\ell}(k)S_{\ell}(k,t)
\label{eq21}%
\end{equation}

In the wavenumber representation, the governing equation of the wave
corresponding to Eq.\ (\ref{eq8}) is%
\begin{equation}
\mathcal{L}_{k}\,S(k,t)=%
\begin{cases}
0 & \text{for }t<-T,\\
F(k,t) & \text{for }-T<t<0,\\
0 & \text{for }t>0.
\end{cases}
\label{eq:40-6}%
\end{equation}
where the wave equation is
\begin{equation}
\mathcal{L}_{k}=e^{-ikx}\,\mathcal{L}\ e^{ikx}%
\end{equation}
and where $F(k,t)$ is the wavenumber representation of the source term
$f(x,t)$
\begin{equation}
F(k,t)\,\,=\,\frac{1}{\sqrt{2\pi}}\int f(x,t)e^{-ikx}\,dx
\end{equation}

Analogous to Eqs.\ \eqref{eq:20230309_105453}, \eqref{eq:20230309_105727}, and
\eqref{eq:20230309_105545}, in the wavenumber representation the wave may be
expressed as the superposition of wavenumber modes
\begin{equation}
S(k,t)\,=\sum_{\ell=1}^{M}S_{\ell}^{B}(k,t),\qquad\text{for }t<-T
\label{eq:20230309_110244}%
\end{equation}
\emph{before} the source is on, a particular solution
\begin{equation}
S(k,t)=\varphi(k,t)\,\qquad\text{for }-T<t<0
\end{equation}
\emph{while} the source is on, where $\varphi$ is related to the spatial
particular solution by
\begin{equation}
\varphi(k,t)=\frac{1}{\sqrt{2\pi}}\int e^{-ikx}\mu(x,t)dx
\label{eq:20231029_130554}%
\end{equation}
and again in terms of modes
\begin{equation}
S(k,t)\,=\sum_{\ell=1}^{M}S_{\ell}^{A}(k,t)\hspace{0.2in}\text{for }t>0
\end{equation}
\emph{after} the source has been turned off.

If we define the $\mathcal{K}$ operator
\begin{equation}
\mathcal{K}=\frac{1}{i}\frac{\partial}{\partial x}, \label{eq:20240917_154231}%
\end{equation}
then the time-evolution of the spatial modes becomes
\begin{equation}
u_{\ell}(x,t)=e^{-i\omega_{\ell}(\mathcal{K})t}u_{\ell}(x,0)
\end{equation}
Thus,
$e^{-i\omega_{\ell}(\mathcal{K})t}$ is the propagator of mode $\ell$ in the spatial representation.

\section{Obtaining the modes from the wave}

\label{sec:obt_modes_frm_wave}

We now show how to obtain the modes when the initial conditions are given by
$u(x,0)$ and its time derivatives. This problem has been
previously addressed \cite{whit,modes}, and we give here a summery. For the
sake of simplicity, we assume that there are two modes; the $M$-mode case is
given in App.\ B. 

Suppose we know $u(x,-T)$ and its time derivatives,
then equivalently we know the initial conditions for the wave in wavenumber space, $S(k,-T)$
and its time derivatives.
For the two-mode case, the
two initial conditions could be expressed in terms of modes. They
are
\begin{equation}
S(k,-T)=S_{1}(k,-T)+S_{2}(k,-T) \label{eq:64}%
\end{equation}
and
\begin{equation}
i\frac{{\partial S(k,-T)}}{{\partial}{t}}=\omega_{1}(k) S_{1}%
(k,-T)+\omega_{2}(k)S_{2}(k,-T) \label{eq:65}%
\end{equation}
where we have used Eq.\ \eqref{eq21}.
Notice that (importantly)
in terms of modes,
no derivatives appear in
the initial conditions
for the wave
(Eqs.\ \eqref{eq:64} and \eqref{eq:65}).
That is,
no derivatives appear on the right-hand side.

Solving the above two equations for the modes, $S_{1}(k,0)$ and $S_{2}(k,0)$, we obtain%
\begin{align}
S_{1}(k,-T)  &  =\frac{\omega_{2}(k)-i\frac{{\partial}}{{\partial}{t}}}%
{\omega_{2}(k)-\omega_{1}(k)}S(k,-T)\label{eq:20230614_162950}\\
S_{2}(k,-T)  &  =\frac{\omega_{1}(k)-i\frac{{\partial}}{{\partial}{t}}}%
{\omega_{1}(k)-\omega_{2}(k)}S(k,-T) \label{eq:20230614_162952}%
\end{align}
The equations may be transformed to the spatial domain to obtain (see
Sec.\ \ref{sec:spatial_formulation})
\begin{align}
u_{1}(x,-T)  &  =\frac{\omega_{2}\left(  \mathcal{K}\right)  -i\frac
{{\partial}{}}{{\partial}{t}}}{\omega_{2}\left(  \mathcal{K}\right)
-\omega_{1}\left(  \mathcal{K}\right)  }u(x,-T)\label{eq:37}\\
u_{2}(x,-T)  &  =\frac{\omega_{1}\left(  \mathcal{K}\right)  -i\frac
{{\partial}{}}{{\partial}{t}}}{\omega_{1}\left(  \mathcal{K}\right)
-\omega_{2}\left(  \mathcal{K}\right)  }u(x,-T) \label{eq:20230614_163034}%
\end{align}

\section{Effect of sources on modes: Wavenumber representation}

As in the previous section, we consider the two-mode case for the sake of
clarity. Also, we use the abbreviated notation $\omega_{1}=\omega_{1}(k)$ and
$\omega_{2}=\omega_{2}(k).$ While the source is on, the wave, $\varphi(k,t),$
is a superposition of the homogeneous part, $S^{B}(k,t)$, and a particular
solution, $S_{p}(k,t)$%

\begin{equation}
\varphi(k,t)=S^{B}(k,t)+S_{p}(k,t)\hspace{0.2in}\hspace{0.2in}\text{for
}-T<t<0 \label{eq32}%
\end{equation}
The solution for all time is then given by Eq.\ \eqref{eq:40-6}, which we now
repeat%
\begin{equation}
S(k,t)=%
\begin{cases}
S^{B}(k,t) & \text{for }t<-T,
\\
\varphi(k,t) & \text{for }-T<t<0,
\\
S^{A}(k,t) & \text{for }t>0.
\end{cases}
\label{126_38}
\end{equation}
where%
\begin{align}
S^{B}(k,t)  &  =S_{1}^{B}(k,t)\,+S_{2}^{B}(k,t)\\
S^{A}(k,t)  &  =S_{1}^{A}(k,t)\,+S_{2}^{A}(k,t)
\end{align}

After the source is turned off, the wave is again expressed as the
superposition of modes. Evolving each mode from that moment by way of
Eq.\ \eqref{eq:20230309_105924}, we obtain that
\begin{equation}
S(k,t)=e^{-i\omega_{1}t}S_{1}^{A}(k,0)+e^{-i\omega_{2}t}S_{2}^{A}%
(k,0)\hspace{0.2in}\text{for }t>0. \label{eq41}%
\end{equation}

At the moment the source is turned off (at time $t=0$) we have
\begin{align}
S_{1}^{A}(k,0)\,+S_{2}^{A}(k,0)\,  &  =\varphi(k,0)\\
\frac{{\partial}{}}{{\partial}{t}}S_{1}^{A}(k,0)\,+\frac{{\partial}{}%
}{{\partial}{t}}S_{2}^{A}(k,0)\,  &  =\frac{{\partial}\varphi(k,0){}%
}{{\partial}{t}}%
\end{align}
and further, from Eq. (\ref{eq21})\ we have that
\begin{equation}
\omega_{1}S_{1}^{A}(k,0)\,+\omega_{2}S_{2}^{A}(k,0)\,=i\frac{{\partial}%
\varphi(k,0){}}{{\partial}{t}}%
\end{equation}
Solving for $S_{2}^{A}(k,0)$ and $S_{2}^{A}(k,0)$ we obtain,
\begin{align}
S_{1}^{A}(k,0)  &  =\frac{\omega_{2}-i\frac{{\partial}{}}{{\partial}{t}}%
}{\omega_{2}-\,\omega_{1}}\varphi(k,0)\label{eq45}\\
S_{2}^{A}(k,0)  &  =\frac{\omega_{1}-i\frac{{\partial}{}}{{\partial}{t}}%
}{\omega_{1}-\,\omega_{2}}\varphi(k,0) \label{eq46}%
\end{align}
Eqs. (\ref{eq45}) and (\ref{eq46}) give the total contribution to modes 1
and 2 by the source term $F(k,t)$ from time $t=-T$ to time $t=0$.

Isolating the change in modes $1$ and $2$ only,
we define
\begin{align}
\Delta S_{1}(k,0) = S_{1}^{A}(k,0) - S_{1}^{B}(k,0) = \frac{ \omega_{2} - i
\frac{ {\partial} } { {\partial}t } } { \omega_{2}-\omega_{1} } S_{p}(k,0)\\
\Delta S_{2}(k,0) = S_{2}^{A}(k,0) - S_{2}^{B}(k,0) = \frac{ \omega_{1} - i
\frac{ {\partial} } { {\partial}t } } { \omega_{1}-\omega_{2} } S_{p}(k,0)
\end{align}


\section{The modes after the source term is turned off}

Now that we know how to obtain the modes from the wave while the source is on, we can evolve the modes from the time it was turned off. Given the
modes $S_{\ell}^{A}(k,0)$ at time 0, we can propagate them forward in time to
time $t>0$. Using Eqs.\ (\ref{eq45}) and (\ref{eq46}),
the modes at some time $t$ after the source has been turned off are
\begin{equation}
S_{1}^{A}(k,t)=e^{-i\omega_{1}t}S_{1}^{A}(k,0)=e^{-i\omega_{1}t}\frac
{\omega_{2}-i\frac{{\partial}{}}{{\partial}{t}}}{\omega_{2}-\,\omega_{1}%
}\varphi(k,0) \label{eq47}%
\end{equation}
\begin{equation}
S_{2}^{A}(k,t)=e^{-i\omega_{2}t}S_{2}^{A}(k,0)=e^{-i\omega_{2}t}\frac
{\omega_{1}-i\frac{\partial}{\partial t}}{\omega_{1}-\omega_{2}}\varphi(k,0)
\label{eq48}%
\end{equation}
or isolating the difference,
\begin{align}
\Delta S_{1}(k,t) = S_{1}^{A}(k,t) - S_{1}^{B}(k,t) = e^{-i\omega_{1}t}
\frac{\omega_{2}-i\frac{{\partial}{}}{{\partial}{t}}} {\omega_{2}-\,\omega_{1}
} S_{p}(k,0)\\
\Delta S_{2}(k,t) = S_{2}^{A}(k,t) - S_{2}^{B}(k,t) = e^{-i\omega_{2}t}
\frac{\omega_{1}-i\frac{{\partial}{}}{{\partial}{t}}} {\omega_{1}-\,\omega_{2}
} S_{p}(k,0)
\end{align}

One can verify that at time $t=0$, the superposition of the two modes is
\begin{equation}
S_{1}(k,0)+S_{2}(k,0)=\varphi(k,0),
\end{equation}
Eqs.\ (\ref{eq47}) and (\ref{eq48}) show that the source-term contribution
to the wave evolves like the homogeneous modes. Meaning that the excitation
due to the source term also evolves like modes do. It has one part that
evolves like mode 1 and one that evolves like mode 2, but only after the
source is turned off.

Thus,
while the source is on,
the particular solution
cannot be expressed in terms of modes,
but after the source is turned off,
it turns into the superposition of two modes.

\subsection{Evolution of the wave after the source term is turned off}

Now that we have shown how to get the modes any time after the source is turned off,
we obtain the wave at $t>0$. For times $t>0$, we have that the superposition of the
modes is
\begin{equation}
S^{A}(k,t)=e^{-i\omega_{1}t}S_{1}^{A}(k,0)+e^{-i\omega_{2}t}S_{2}%
^{A}(k,0)\hspace{0.2in}\text{for }t>0
\end{equation}
Using Eqs.\ Eq. (\ref{eq45}) and Eq. (\ref{eq46}) we have that explicitly
\begin{align}
S^{A}(k,t)  &  =\left(  e^{-i\omega_{1}t}\frac{\omega_{2}-i\frac{{\partial}{}%
}{{\partial}{t}}}{\omega_{2}-\,\omega_{1}}+e^{-i\omega_{2}t}\frac{\omega
_{1}-i\frac{{\partial}{}}{{\partial}{t}}}{\omega_{1}-\,\omega_{2}}\right)
\varphi(k,0)\\
&  =-e^{-i(\omega_{1}+\omega_{2})t}\frac{\partial}{\partial t}e^{i\frac
{\omega_{1}+\omega_{2}}{2}t}\,\frac{\sin\left(  \frac{\omega_{1}-\omega_{2}%
}{2}t\right)  }{\left(  \frac{\omega_{1}-\omega_{2}}{2}\right)  }%
\varphi(k,t=0) \label{eq:20231026_132640}%
\end{align}

\section{Spatial Formulation}

\label{sec:spatial_formulation}

In the spatial domain, the wave corresponding to Eq.\ (\ref{126_38}) is
\begin{equation}
u(x,t)=%
\begin{cases}
u^{B}(x,t) & \text{for }t<-T,\\
\mu(x,t) & \text{for }-T<t<0,\\
u^{A}(x,t) & \text{for }0<t.
\end{cases}
\end{equation}
and as in Eq.\ \eqref{eq32},
\begin{align}
\mu(x,t)
=
u^B(x,t)
+
u_p(x,t)
\text{, for }-T<t<0
%
\end  {align}

Eqs.\ (\ref{eq47})\ and (\ref{eq48})\ relate the wavenumber-domain modes after
the source has been turned off to those before it was on. To obtain the
equivalent relation in the spatial domain, we Fourier transform
Eq.\ (\ref{eq48}) to obtain that
\begin{equation}
u_{1}^{A}(x,t)=\frac{1}{2\pi}\iint e^{-i\omega_{1}t}e^{ik(x-x^{\prime})}%
\frac{\omega_{2}(k)-i\frac{{\partial}}{{\partial}{t}}}{\omega_{2}%
(k)-\,\omega_{1}(k)}\mu(x^{\prime},0)dx^{\prime}dk \label{eq72}%
\end{equation}
for mode 1. Similarly, for mode 2, Eq.\ (\ref{eq48}) gives
\begin{equation}
u_{2}^{A}(x,t)=\frac{1}{2\pi}\iint e^{-i\omega_{2}t}e^{ik(x-x^{\prime})}%
\frac{\omega_{1}(k)-i\frac{{\partial}}{{\partial}{t}}}{\omega_{1}%
(k)-\,\omega_{2}(k)}\mu(x^{\prime},0)dx^{\prime}dk \label{eq74}%
\end{equation}
Alternatively, using the differential operator $\mathcal{K}$, 
modes 1 and 2 in the spatial domain are
\begin{equation}
u_{1}^{A}(x,t)=e^{-i\omega_{1}(\mathcal{K})t}\frac{\omega_{2}(\mathcal{K}%
)-i\frac{\partial}{\partial t}}{\omega_{2}(\mathcal{K})-\omega_{1}%
(\mathcal{K})}\mu(x,0) \label{eq:20231026_131302}%
\end{equation}
and
\begin{equation}
u_{2}^{A}(x,t)=e^{-i\omega_{2}(\mathcal{K})t}\frac{\omega_{1}(\mathcal{K}%
)-i\frac{\partial}{\partial t}}{\omega_{1}(\mathcal{K})-\omega_{2}%
(\mathcal{K})}\mu(x,0) \label{eq:20231026_131240}%
\end{equation}

The spatial wave at times $t>0$ is then the superposition of the modes
\begin{equation}
u^{A}(x,t)=u_{1}^{A}(x,t)+u_{2}^{A}(x,t)\hspace{0.2in}\text{for }t>0
\end{equation}
and using Eqs.\ (\ref{eq74}) and (\ref{eq72}), we obtain
\begin{equation}
u^{A}(x,t)=\frac{i}{2\pi}\iint\frac{e^{-i[\omega_{1}(k)+\omega_{2}(k)]t}%
}{\omega_{1}(k)-\omega_{2}(k)}e^{ik(x-x^{\prime})}\frac{\partial}{\partial
t}\left(  e^{i\omega_{2}(k)t}-e^{i\omega_{1}(k)t}\right)  \mu(x^{\prime
},t=0)\,dx^{\prime}dk
\end{equation}
Alternatively using Eqs.\ \eqref{eq:20231026_131302} and
\eqref{eq:20231026_131240} we also have%
\begin{equation}
u^{A}(x,t)=i\frac{e^{-i[\omega_{1}(\mathcal{K})+\omega_{2}(\mathcal{K})]t}%
}{\omega_{1}(\mathcal{K})-\omega_{2}(\mathcal{K})}\frac{\partial}{\partial
t}\left(  e^{i\omega_{2}(\mathcal{K})t}-e^{i\omega_{1}(\mathcal{K})t}\right)
\mu(x,t=0)
\end{equation}

\section{ Special case: Source-free case}

Treating the source-free case in general, one takes
\begin{equation}
u_{p}(x,t)=0\qquad\text{for }-T<t<0
\end{equation}
Then when the source is turned off at $t=0$,
\begin{equation}
\mu(x,0)=u^{A}(x,0) \label{126_62}%
\end{equation}
and during $-T \le t \le0$
\begin{equation}
u^{B}(x,t)=\mu(x,t)
\end{equation}
Since in the source-free case, we are allowed to set the source duration, $T$,
to be zero, we can identify time $t=-T$ with time $t=0$, and we thus also
have
\begin{equation}
u^{B}(x,0)=\mu(x,0) \label{126_64}%
\end{equation}
Therefore, from Eqs.\ \eqref{126_62} and \eqref{126_64} we have
\begin{equation}
u^{B}(x,0)=u^{A}(x,0)
\end{equation}
and as a result,
\begin{align}
u^{A}(x,t)  &  =u^{B}(x,t)\\
&  =u_{1}^{B}(x,t)+u_{2}^{B}(x,t)\\
&  =e^{-i\omega_{1}(\mathcal{K})t}u_{1}^{B}(x,0)+e^{-i\omega_{2}%
(\mathcal{K})t}u_{2}^{B}(x,0).
\end{align}

\section{Special case: Single mode wave equation}

There are cases where the wave equation has only one mode. One example is the
free particle Schr\"{o}dinger-type equation
\begin{equation}
\frac{\partial u}{\partial t}=i\gamma\frac{\partial^{2}u}{\partial x^{2}}%
\end{equation}
where $\gamma=\hbar/2m$, in which case the dispersion relation is
$\omega(k)=\gamma k^{2}.$ Another example is the linearized Korteweg-de Vries
equation {\cite{whit}}
\begin{equation}
\nu\frac{\partial^{3}u}{\partial x^{3}}+c_{0}\frac{\partial u}{\partial
x}=-\frac{\partial u}{\partial t}%
\end{equation}
which gives the single mode $\omega(k)=c_{0}k-\nu k^{3}$.

All of these wave equations have the following solution.
After time $t=0$,
the solution is
\begin{equation}
S^{A}(k,t)=e^{-i\omega(k)t}\varphi(k,0),
\end{equation}
and in the position representation,
\begin{equation}
u^{A}(x,t)=\frac{1}{2\pi}%
{\displaystyle\iint}
e^{-i\omega(k)t}e^{ik(x-x^{\prime})}\mu(x^{\prime},0)dkdx^{\prime}%
\end{equation}
or in terms of the $\mathcal{K}$ operator,
we have
\begin{equation}
u^{A}(x,t)=e^{-i\omega(\mathcal{K})t}\mu(x,0)
\end{equation}

\section{Special case: Two-mode case where $\omega_{1}(k)=-\omega_{2}(k)$}

It is common in the two-mode case that the two solutions of the dispersion
equation, Eq.\ \eqref{eq03}, are related by
\begin{equation}
\omega_{1}(k)=-\omega_{2}(k) \label{eq:20230314_135722}%
\end{equation}
Among such equations are the standard wave equation, the beam equation, and
the linear Boussinesq equation {\cite{whit}}. Specifically for the standard
wave equation,
\begin{equation}
\frac{\partial^{2}u}{\partial x^{2}}=\frac{1}{c^{2}}\frac{\partial^{2}%
u}{\partial t^{2}},
\end{equation}
the modes are
\begin{equation}
\omega_{1}(k)=-\omega_{2}(k)=ck
\end{equation}
The standard wave equation is discussed in Sec.\ \ref{sec:wave_eqn}.

For the beam equation \cite{whit},
\begin{equation}
\frac{\partial^{4}u}{\partial x^{4}}=-\frac{1}{\gamma^{2}}\frac{\partial^{2}%
u}{\partial t^{2}}%
\end{equation}
the modes are
\begin{equation}
\omega_{1}(k)=-\omega_{2}(k)=\gamma^{2}k^{2},
\end{equation}
and for the linear Boussinesq equation \cite{whit},
\begin{equation}
\frac{\partial^{2}u}{\partial t^{2}}-c^{2}\frac{\partial^{2}u}{\partial x^{2}%
}-\beta^{2}\frac{\partial^{4}u}{\partial x^{2}\partial t^{2}}=0
\end{equation}
the modes are given by
\begin{equation}
\omega_{1}(k)=-\omega_{2}(k)=\frac{c^{2}k}{\sqrt{1+\beta^{2}k^{2}}}.
\end{equation}

\subsection{Individual wavenumber modes}


We now apply the results of the previous section to the situation $\omega
_{1}(k)=-\omega_{2}(k)$, where we use the abbreviated notation $\omega = \omega_1(k)$. 
We start by considering each of the modes in wavenumber space
individually. Using Eqs.\ (\ref{eq47}) and (\ref{eq48}) we have%

\begin{equation}
S_{1}^{A}(k,t)=\frac{1}{2}e^{-i\omega t}\left[  1-\frac{1}{i\omega}%
\frac{{\partial}{}}{{\partial}{t}}\right]  \varphi(k,0) \label{eq86}%
\end{equation}
and
\begin{equation}
S_{2}^{A}(k,t)=\frac{1}{2}e^{i\omega t}\left[  1+\frac{1}{i\omega}%
\frac{{\partial}{}}{{\partial}{t}}\right]  \varphi(k,0) \label{eq87}%
\end{equation}

\subsection{The wave in wavenumber space}

We find that in wavenumber space, the wave is the superposition of
Eqs.\ \eqref{eq86} and \eqref{eq87}. It could also be obtained by plugging in
$\omega=\omega_{1}=-\omega_{2}$
in Eq.\ \eqref{eq:20231026_132640},%

\begin{align}
S^{A}(k,t)  &  =\left(  \cos\omega t\,+\frac{\sin\omega t}{\omega}%
\frac{{\partial}{}}{{\partial}{t}}\right)  \varphi(k,0)\\
&  =\frac{\partial}{\partial t}\frac{\sin\,\omega t}{\omega}\varphi(k,t=0)
\end{align}
or explicitly in terms of $S_{p}$,
%
%
\begin{align}
\Delta S(k,t)  &  =\left(  \cos\omega t\,+\frac{\sin\omega t}{\omega}%
\frac{{\partial}{}}{{\partial}{t}}\right)  S_{p}(k,0)\\
&  =\frac{\partial}{\partial t}\frac{\sin\,\omega t}{\omega}S_{p}(k,t=0)
\end{align}

\subsection{Position space: Modes}


In position space, using Eq.\ \eqref{eq74} the modes are
\begin{equation}
u_{1}^{A}(x,t)=\frac{1}{4\pi}\iint e^{-i\omega(k)t}e^{ik(x-x^{\prime})}\left(
1+\frac{i}{\omega(k)}\frac{\partial}{\partial t}\right)  \mu(x^{\prime
},0)dkdx^{\prime}%
\end{equation}
or using Eq.\ \eqref{eq:20231026_131302}
\begin{equation}
u_{1}^{A}(x,t)=\frac{-1}{i\omega(\mathcal{K})}\frac{\partial}{\partial
t}e^{-i\omega(\mathcal{K})t}\mu(x,t=0)
\end{equation}
Similarly,
for mode 2,
\begin{equation}
u_{2}^{A}(x,t)=\frac{1}{4\pi}\iint e^{i\omega(k)t}e^{ik(x-x^{\prime})}\left(
1-\frac{i}{\omega(k)}\frac{\partial}{\partial t}\right)  \mu(x^{\prime
},0)dkdx^{\prime}%
\end{equation}
or
\begin{equation}
u_{2}^{A}(x,t)=\frac{1}{i\omega(\mathcal{K})}\frac{\partial}{\partial
t}e^{i\omega(\mathcal{K})t}\mu(x,t=0)
\end{equation}

\subsection{Position space: The wave from modes in position space}

The wave in position space is
\begin{align}
u^{A}(x,t)  &  =\frac{1}{2\pi}\iint e^{ik(x-x^{\prime})}\left(  \cos\left(
\omega(k)t\right)  +\frac{\sin\left(  \omega(k)t\right)  }{\omega(k)}%
\frac{{\partial}}{{\partial}{t}}\right)  \mu(x^{\prime},0)dkdx^{\prime
}\label{eq:45-67}\\
&  =\frac{1}{2\pi}\iint e^{ik(x-x^{\prime})}\frac{\partial}{\partial t}%
\frac{\sin\,\omega(k)t}{\omega(k)}\mu(x^{\prime},t=0)dkdx^{\prime}
\label{eq:45-67}%
\end{align}
or
\begin{align}
u^{A}(x,t)  &  =\frac{-1}{i\omega(\mathcal{K})}\frac{\partial}{\partial
t}e^{-i\omega(\mathcal{K})t}\mu(x,t=0)+\frac{1}{i\omega(\mathcal{K})}%
\frac{\partial}{\partial t}e^{i\omega(\mathcal{K})t}\mu(x,t=0)\\
&  =\frac{1}{i\omega(\mathcal{K})}\frac{\partial}{\partial t}\left[
e^{i\omega(\mathcal{K})t}-e^{-i\omega(\mathcal{K})t}\right]  \mu(x,t=0)
\end{align}
giving that
\begin{align}
u^{A}(x,t) 
=2\frac{\partial}{\partial t}\frac{\sin\left[  \omega(\mathcal{K})t\right]
}{\omega(\mathcal{K})}\mu(x,t=0)
\end{align}

\section{Generalization of the d'Alembert solution}

Throughout the paper,
we have investigated the case
of a finite-duration source;
we found that this framework clarifies many issues.
However,
for providing the generalization of the d'Alembert solution,
(which we do in this section)
we decided to deviate from this framework
to allow our results could be more readily comparable
with the usual case considered by the d'Alembert solution.
The usual d'Alembert case
is a special case of our case,
which we present
in App.\ A. 

%

In this section,
we provide a derivation of the d'Alembert solution
for an arbitrary two-mode wave equation.
For the standard wave equation, the solution for the wave may be explicitly
expressed in terms of the wave and its time derivative at time zero \cite{Davis,Bland}.
The solution is called the d'Alembert solution. The standard wave equation has two
modes which are negative of each other, and its dispersion relation is linear
in $k$. We now consider the case where the two modes are arbitrary.

In App.\ A, 
we show that for an arbitrary two-mode wave equation, the solution in wavenumber space is
\begin{align}
S(k,t)  &  =\frac{e^{-i\omega_{1}t}\omega_{2}(k)-e^{-i\omega_{2}t}\omega
_{1}(k)}{\omega_{2}(k)-\omega_{1}(k)}S(k,0)+i\frac{\left(  -e^{-i\omega_{1}%
t}+e^{-i\omega_{2}t}\right)  }{\omega_{2}(k)-\omega_{1}(k)}\frac{{\partial}%
}{{\partial}{t}}S(k,0)\nonumber\\
&  -i\int_{0}^{t}\left(  -\frac{e^{-i\omega_{1}(t-t^{\prime})}}{\omega
_{2}(k)-\omega_{1}(k)}+\frac{e^{-i\omega_{2}(t-t^{\prime})}}{\omega
_{2}(k)-\omega_{1}(k)}\right)  F(k,t^{\prime})dt^{\prime} \label{126_145_2}%
\end{align}
and in position space it is
\begin{align}
u(x,t)  &  =\frac{1}{2\pi}%
{\displaystyle\iint}
\left(  \frac{e^{-i\omega_{1}t-ik(x^{\prime}-x)}\omega_{2}(k)}{\omega
_{2}(k)-\omega_{1}(k)}-\frac{e^{-i\omega_{2}t-ik(x^{\prime}-x)}\omega_{1}%
(k)}{\omega_{2}(k)-\omega_{1}(k)}\right)  u(x^{\prime},0)dkdx^{\prime
}\nonumber\\
&  -\frac{1}{2\pi i}\iint \left(  -\frac{e^{-i\omega_{1}t}}{\omega
_{2}(k)-\omega_{1}(k)}+\frac{e^{-i\omega_{2}t}}{\omega_{2}(k)-\omega_{1}%
(k)}\right)  e^{-ik(x^{\prime}-x)}\frac{{\partial}%
}{{\partial}{t}}u(x^{\prime},0)dkdx^{\prime}\nonumber\\
&  +\frac{1}{2\pi i}\int_{0}^{t}\left[
{\displaystyle\iint}
\left(  -\frac{e^{-i\omega_{1}(t-t^{\prime})}}{\omega_{2}(k)-\omega_{1}%
(k)}+\frac{e^{-i\omega_{2}(t-t^{\prime})}}{\omega_{2}(k)-\omega_{1}%
(k)}\right)  f(x^{\prime},t^{\prime})e^{-ik(x^{\prime}-x)}dx^{\prime
}dk\right]  dt^{\prime} \label{126_146_2}%
\end{align}
Eq.\ \eqref{126_145_2} is the generalization of the d'Alembert solution in wavenumber space. Similarly, Eq.\ \eqref{126_146_2} is the d'Alembert solution in position space.
For any arbitrary two-mode wave equation,
these give the wave at time $t$,
and they reduce to the regular d'Alembert solution
for the case of the standard wave equation;
that is, where the two modes are
$\omega_{1,2}=\pm ck$.
This is shown in the following subsections.
%


\subsection{D'Alembert-type solution for $\omega_{2}(k)=-\omega_{1}(k)$}

A common case is where the two modes are the negative of each other. We take
\begin{equation}
\omega=\omega_{1}(k)=-\omega_{2}(k)
\end{equation}
Specializing the above solutions for this case we obtain in wavenumber space
that%
\begin{align}
S(k,t)  &  =\frac{1}{2}\left(  e^{-i\omega t}+e^{i\omega t}\right)
S(k,0)+\frac{i}{2\omega (k)}\left(  e^{-i\omega t}+e^{i\omega %
t}\right)  \frac{{\partial}}{{\partial}{t}}S(k,0)\nonumber\\
&  -i\frac{1}{2}\int_{0}^{t}\left(  \frac{e^{-i\omega (t-t^{\prime}%
)}-e^{i\omega (t-t^{\prime})}}{\omega (k)}\right)  F(k,t^{\prime
})dt^{\prime} \\ 
&  =-i\sin(\omega t)S(k,0)+i\frac{\cos(\omega t)}{\omega}\frac{\partial
}{\partial t}S(k,0)-\int\frac{\sin[\omega(t-t^{\prime})]}{\omega}%
F(k,t^{\prime})dt^{\prime} \label{126_148}%
\end{align}
and in position space we have%

\begin{align}
u(x,t)  &  =\frac{1}{2\pi}\frac{1}{2}%
{\displaystyle\iint}
e^{-ik(x^{\prime}-x)}\left(  e^{-i\omega t}+e^{i\omega t}\right)
u(x^{\prime},0)dkdx^{\prime}\nonumber\\
&  -\frac{1}{2\pi i}\frac{1}{2}\iint\frac{\left(  e^{-i\omega t}%
+e^{i\omega t}\right)  }{\omega (k)}e^{-ik(x^{\prime}-x)} \frac{{\partial}}{{\partial}{t}}u(x^{\prime}%
,0)dkdx^{\prime}\nonumber\\
&  +\frac{1}{2\pi i}\frac{1}{2}\int_{0}^{t}\left[
{\displaystyle\iint}
\left(  \frac{e^{-i\omega (t-t^{\prime})}-e^{i\omega (t-t^{\prime})}%
}{\omega (k)}\right)  f(x^{\prime},t^{\prime})e^{-ik(x^{\prime}%
-x)}dx^{\prime}dk\right]  dt^{\prime}%
\label{eq:20250324_132619}
\end{align}
\bigskip

\subsection{Standard wave equation}
\label{sec:wave_eqn}

It is of interest to show how the above equations reduce to the standard
d'Alembert solution. The standard wave equation with a source term is
\begin{equation}
\frac{\partial^{2}u}{\partial x^{2}}-\frac{1}{c^{2}}\frac{\partial^{2}%
u}{\partial t^{2}}=%
\begin{cases}
0 & \text{for }t<-T,\\
f(x,t) & \text{for }-T<t<0,\\
0 & \text{for }t>0.
\end{cases}
\label{126_97}%
\end{equation}
In the wavenumber representation it is
\begin{equation}
(ik)^{2}S-\frac{1}{c^{2}}\frac{\partial^{2}S}{\partial t^{2}}=%
\begin{cases}
0 & \text{for }t<-T,\\
F(k,t) & \text{for }-T<t<0,\\
0 & \text{for }t>0.
\end{cases}
\end{equation}
The dispersion relation is
\begin{equation}
\omega^{2}=c^{2}k^{2}%
\end{equation}
and the two modes are given by
\begin{equation}
\omega_{1}(k)=ck\hspace{0.5in}\omega_{2}(k)=-ck \label{eq106}%
\end{equation}
Substituting these values in Eq.\ (\ref{126_148}), we obtain
\begin{align}
S(k,t)  &  =\frac{1}{2}\left(  e^{-ickt}+e^{ickt}\right)  S(k,0)+\frac{i}%
{2ck}\left(  e^{-ickt}-e^{ickt}\right)  \frac{{\partial}}{{\partial}{t}%
}S(k,0)\nonumber\\
&  +\frac{1}{2ic}\int_{0}^{t}\left(  \frac{e^{-ick(t-t^{\prime})}%
-e^{ick(t-t^{\prime})}}{k}\right)  F(k,t^{\prime})dt^{\prime}\\
&  =\cos(ckt)S(k,0)+\frac{\sin(ckt)}{ck}\frac{\partial}{\partial t}%
S(k,0)-\int_{0}^{t}\frac{\sin[ck(t-t^{\prime})}{ck}F(k,t^{\prime})dt^{\prime}%
\end{align}
In position space, we use the modes $\omega$ in Eq.\ \eqref{eq:20250324_132619}
(and integrate by parts)
\begin{align}
u(x,t)  &  =\frac{1}{2\pi}\frac{1}{2}\int\left(  e^{-ickt}+e^{+ickt}\right)
e^{-ik(x^{\prime}-x)}
u(x^{\prime},0)dkdx^{\prime
}\nonumber\\
&  -\frac{1}{2\pi}\frac{1}{2c}\iint\left(  e^{-ickt}-e^{ickt}\right)
e^{-ik(x^{\prime}-x)}\int_{\alpha}^{x^{\prime}}\frac{{\partial}}{{\partial}%
{t}}u(x^{\prime\prime},0)dx^{\prime\prime}dkdx^{\prime}\nonumber\\
&  +\frac{1}{2\pi}\frac{1}{2c}\int_{0}^{t}\left[
{\displaystyle\iint}
\left(  e^{-ick(t-t^{\prime})}-e^{ick(t-t^{\prime})}\right)
e^{-ik(x^{\prime}-x)}  \int_{\alpha}^{x^{\prime}}f(x^{\prime\prime
},t^{\prime})dx^{\prime\prime}dx^{\prime}dk\right]  dt^{\prime}%
\end{align}
which simplifies to
\begin{align}
u(x,t)  &  =\frac{1}{2}\int\left( \Big.  \delta(x^{\prime}-x+ct)+\delta(x^{\prime
}-x-ct)\right)  u(x^{\prime},0)dx^{\prime}\nonumber\\
&  -\frac{1}{2c}\int\left(  \Big. \delta(x^{\prime}-x+ct)-\delta(x^{\prime
}-x-ct)\right)  \int_{\alpha}^{x^{\prime}}\frac{{\partial}}{{\partial}{t}%
}u(x^{\prime\prime},0)dkdx^{\prime}\nonumber\\
&  +\frac{1}{2c}\int_{0}^{t}\left[
{\displaystyle\iint}
\left(  \Big. \delta(x^{\prime}-x+ck(t-t^{\prime}))-\delta(x^{\prime}%
-x-ck(t-t^{\prime}))\right)  \int_{\alpha}^{x^{\prime}}f(x^{\prime\prime
},t^{\prime})dx^{\prime\prime}dx^{\prime}dk\right]  dt^{\prime}%
\end{align}
giving
\begin{align}
u(x,t)  &  =\frac{1}{2}\left(  u(x+ct,0)+u(x-ct,0)\Big.\right)  +\frac{1}%
{2c}\int_{x-ct}^{x+ct}\frac{{\partial}}{{\partial}{t}}u(x^{\prime\prime
},0)\nonumber\\
&  -\frac{1}{2c}\int_{0}^{t}\int_{x-c(t-t^{\prime})}^{x+c(t-t^{\prime}%
)}f(x^{\prime\prime},t^{\prime})dx^{\prime\prime}dt^{\prime}%
\label{eq:20250324_132541}
\end{align}
which is the well-known d'Alembert solution for the case of the standard wave
equation with sources \cite{pinch}.

\bigskip

\section{Conclusion and Summary}

We have obtained explicit equations for the change in the modes when a
finite-duration source is applied.
%
We have also generalized the d'Alembert solution to wave equations other than
the standard wave equation.
In particular to equations with arbitrary number of modes and arbitrary
relations between them.

We now summarize the results for the two-mode case.
The $M$-mode case is given in App.\ B.
%
The source term is turned on at time $t=-T$, and turned off at time $t=0$.
Before the source is turned on, the wave is given by the superposition of
modes, $u_{\ell}^{B}$,
\begin{equation}
u(x,t)\,=\sum_{\ell=1}^{M}u_{\ell}^{B}(x,t),\qquad\text{for }t<-T
\end{equation}
and after the source by modes given by $u_{\ell}^{A}(x,t)$
\begin{equation}
u(x,t)\,=\sum_{\ell=1}^{M}u_{\ell}^{A}(x,t),\hspace{0.2in}\text{for }t>0
\end{equation}
While the source is on, that is, between times $t=-T$ to $t=0$, the wave, $u(x,t)$, is
given by
\begin{equation}
u(x,t)=\mu(x,t)\qquad\text{for }-T<t<0
\end{equation}
where $\mu(x,t)$ is the solution to the wave equations during that time.

In wavenumber space, the solution, $S^{A}(k,t)$, after the source is turned
off is
\begin{equation}
S^{A}(k,t)=\sum_{\ell=1}^{M}e^{-i\omega_{\ell}t}S_{\ell}^{A}(k,0)
\end{equation}
where the modes $S_{\ell}(k,0)$ are given by Eq.\ \eqref{eq:20240917_153315}.

In position space we have
\begin{equation}
u^A(x,t)=\sum_{\ell=1}^{M}e^{-i\omega_{\ell}(\mathcal{K})t}u_{\ell}^{A}(x,0)
\end{equation}
where the $u_{\ell}^{A}(x,0)$ are given by Eq.\ \eqref{eq:20240917_153426},
and the $\mathcal{K}$ operator in Eq.\ \eqref{eq:20240917_154231}. These modes
contain the effect of the source term. We have applied the method to a number
of special cases, including the beam equation, and the linear Boussinesq
equation.

In addition,
using the methods presented,
we obtained a generalization of the d'Alembert solution which applies to an arbitrary linear wave equation.

\appendix

\renewcommand{\thesection}{Appendix \Alph{section}:}

\section{Generalization of the d'Alembert solution}
\label{sec:20241211_144359}

The wave, $u(x,t)$ in position space, and in wavenumber space, is the superposition
of the homogeneous and particular solutions. That is, given the wave at time
$t^{\prime}$ before the source was on, then at time $t$, after the source has
been turned off, $u^{A}(x,t)$ is the superposition of $u^{B}(x,t^{\prime})$
propagated to time $t$, and the change in the modes $\Delta u$.

We treat source free case separately from the source case, and then
superimpose them to find the full solution

\bigskip

\noindent\textbf{Source free case}. In
Sec.\ \ref{sec:obt_modes_frm_wave}, we showed that in wavenumber space the two
modes are given by%
\begin{align}
S_{1}(k,-T)  &  =\frac{\omega_{2}(k)-i\frac{{\partial}}{{\partial}{t}}}%
{\omega_{2}(k)-\omega_{1}(k)}S(k,t=-T)\\
S_{2}(k,-T)  &  =-\frac{\omega_{1}(k)-i\frac{{\partial}}{{\partial}{t}}}%
{\omega_{2}(k)-\omega_{1}(k)}S(k,t=-T)
\end{align}
where $S(k,0)$ and $\frac{{\partial}}{{\partial}{t}}S(k,0)$ are the initial
conditions in wavenumber space. To propagate $S(k,t)$ forward in time, we propagate
the modes
\begin{equation}
S(k,t)=e^{-i\omega_{1}(t+T)}S_{1}(k,-T)+e^{-i\omega_{2}(t+T)}S_{2}(k,-T)
\end{equation}
and therefore we have that the general solution is given by
\begin{equation}
S(k,t)=\frac{e^{-i\omega_{1}(t+T)}\omega_{2}(k)-e^{-i\omega_{2}(t+T)}\omega_{1}%
(k)}{\omega_{2}(k)-\omega_{1}(k)}S(k,-T)-i\frac{e^{-i\omega_{1}(t+T)}%
-e^{-i\omega_{2}(t+T)}}{\omega_{2}(k)-\omega_{1}(k)}\frac{{\partial}}{{\partial
}{t}}S(k,t=-T)
\label{eq135}
\end{equation}
This applies to arbitrary $\omega_{2}(k)$ and $\omega_{1}(k)$. By the way, it
could also be written as
\begin{align}
S(k,t) = -i e^{-i(\omega_{2}+\omega_{1})(t+T)} \frac{{\partial}}{{\partial}{t}}
\left(  \frac{ e^{ i\omega_{2}(t+T)} - e^{ i\omega_{1}(t+T)} } {\omega_{2}%
(k)-\omega_{1}(k)} \right)  S(k,t=-T) \label{eq:20241121_135146}%
\end{align}

We note that while a wave with more than one mode does not have a legitimate propagator,
in some sense, we have succeeded to find a propagator for the full wave,
\begin{align}
K_{k}(k;t,\bar t) = -i e^{-i(\omega_{2}+\omega_{1})(t-\bar t)} \frac
{{\partial}}{{\partial}{t}} \left(  \frac{ e^{ i\omega_{2}(t-\bar t)} - e^{
i\omega_{1}(t-\bar t)} } {\omega_{2}(k)-\omega_{1}(k)} \right)  ,
\end{align}
and the reason that it works is that it contains a derivative (while a legitimate propagator does not contain derivatives).
Here the $t$-derivative is applied to both the terms in the parentheses and to any function to the right;
thus
\begin{align}
S(k,t) = K_{k}(k;t,\bar t) S(k,\bar t)
\end{align}

For the solution in position space, Fourier transform
Eq.\ (\ref{eq135})
to
obtain
\begin{align}
\frac{1}{\sqrt{2\pi}}\int S(k,t)e^{ikx}\,dk  &  =\frac{1}{\sqrt{2\pi}}\int
e^{ikx} \left(  \frac{ e^{-i\omega_{1}(t+T)} \omega_{2}(k) - e^{-i\omega_{2}(t+T)}
\omega_{1}(k) } {\omega_{2}(k)-\omega_{1}(k)} \right)  S(k,-T)dk\nonumber\\
&  -i\frac{1}{\sqrt{2\pi}}\int e^{ikx} \left(  \frac{ e^{-i\omega_{1}(t+T)} -
e^{-i\omega_{2}(t+T)} } {\omega_{2}(k)-\omega_{1}(k)} \right)  \frac{{\partial}%
}{{\partial}{t}}S(k,t=-T)dk
\end{align}
which could also be expressed as
\begin{align}
\int\frac{ e^{ikx} }{\sqrt{2\pi}} S(k,t)\,dk  &  = -i \int\frac{ e^{ikx}
}{\sqrt{2\pi}} e^{-i(\omega_{2}+\omega_{1})(t+T)} \frac{{\partial}}{{\partial}{t}}
\left(  \frac{ e^{ i\omega_{2}(t+T)} - e^{ i\omega_{1}(t+T)} } {\omega_{2}%
(k)-\omega_{1}(k)} \right)  S(k,t=-T)dk
\label{eq:20250211_142603}
\end{align}
and substituting the following expression for $S(k,t)$ in Eq.\ \eqref{eq:20250211_142603}
\begin{equation}
S(k,t=-T)=\frac{1}{\sqrt{2\pi}}\int u(x^{\prime},t=-T)e^{-ikx^{\prime}}dx^{\prime}%
\end{equation}
we obtain
\begin{align}
u(x,t)  &  = {\displaystyle\iint} \frac{ e^{-ik(x^{\prime}-x)} }{2\pi}\left(
\frac{ e^{-i\omega_{1}(t+T)} \omega_{2}(k) - e^{-i\omega_{2}(t+T)} \omega_{1}(k) }
{\omega_{2}(k)-\omega_{1}(k)} \right)  u(x^{\prime},t=-T)dkdx^{\prime}\nonumber\\
&  +i {\displaystyle\iint} \frac{ e^{-ik(x^{\prime}-x)} }{2\pi}\left(  \frac{
e^{-i\omega_{1}(t+T)} - e^{-i\omega_{2}(t+T)} } {\omega_{2}(k)-\omega_{1}(k)} \right)
\frac{{\partial}}{{\partial}{t}}u(x^{\prime},t=-T)dkdx^{\prime}%
\label{eq:20241121_140650}
\end{align}

Thus far, we have only considered the source-free case.
We now obtain the source contribution;
the wave is the superposition of the two.

\bigskip

\noindent\textbf{Source term contribution}. It could
be shown that for any two-mode wave equation, the source term contribution (in
wavenumber space) is
\begin{align}
\Delta S(t)  &  = i \int_{-T}^{0} \frac{ e^{-i\omega_{1}(t-t^{\prime})} -
e^{-i\omega_{2}(t-t^{\prime})} }{\omega_{2}(k)-\omega_{1}(k)} F(k,t^{\prime
})dt^{\prime} \label{126_141}%
\end{align}
Where as usual, the source term in the wavenumber representation,
$F(k,t^{\prime})$, is given by
\begin{equation}
F(k,t^{\prime})=\frac{1}{\sqrt{2\pi}}\int f(x^{\prime},t^{\prime}%
)e^{-ikx}dx^{\prime}%
\end{equation}
Taking the Fourier transform of both sides of Eq.\ (\ref{126_141}), we find
that in position space, the source-term contribution to the wave is
\begin{equation}
\Delta u(t) = \int\frac{e^{ikx}}{\sqrt{2\pi}}\Delta Sdk = -\int_{-T}^{0}
\left[  \iint\frac{ e^{ik(x-x^{\prime})} }{2\pi i} \frac{ e^{-i\omega
_{1}(t-t^{\prime})} - e^{-i\omega_{2}(t-t^{\prime})} } {\omega_{2}%
(k)-\omega_{1}(k)}f(x^{\prime},t^{\prime}) dx^{\prime} dk \right]  dt^{\prime}
\label{eq:20241121_140558}%
\end{equation}

\bigskip

\noindent\textbf{Full solution}. The wave is the
superposition of the source-term contribution and the free propagation part;
in position space it is the superposition of Eqs.\ \eqref{eq:20241121_140650}
and \eqref{eq:20241121_140558}, and in wavenumber space it is the
superposition of Eqs.\ \eqref{eq:20241121_135146} and \eqref{126_141}.
Explicitly, in wavenumber space, the wave is
\begin{align}
S(k,t)  &  =\frac{e^{-i\omega_{1}(t+T)}\omega_{2}(k)-e^{-i\omega_{2}(t+T)}\omega
_{1}(k)}{\omega_{2}(k)-\omega_{1}(k)}S(k,-T)-i\frac{ e^{-i\omega_{1}%
(t+T)}-e^{-i\omega_{2}(t+T)} }{\omega_{2}(k)-\omega_{1}(k)}\frac{{\partial}}%
{{\partial}{t}}S(k,-T)\nonumber\\
&  +i\int_{-T}^{0} \frac{ e^{-i\omega_{1}(t-t^{\prime})} - e^{-i\omega
_{2}(t-t^{\prime})} } {\omega_{2}(k)-\omega_{1}(k)} F(k,t^{\prime})dt^{\prime}
\label{126_145}%
\end{align}
which could also be written as (using Eq.\ \eqref{eq:20241121_135146})
\begin{align}
S(k,t)  &  = -i e^{-i(\omega_{2}+\omega_{1})(t+T)} \frac{{\partial}}{{\partial}%
{t}} \left(  \frac{ e^{ i\omega_{2}(t+T)} - e^{ i\omega_{1}(t+T)} } {\omega
_{2}(k)-\omega_{1}(k)} \right)  S(k,-T)dk\nonumber\\
&  +i\int_{-T}^{0} \frac{ e^{-i\omega_{1}(t-t^{\prime})} - e^{-i\omega
_{2}(t-t^{\prime})} } {\omega_{2}(k)-\omega_{1}(k)} F(k,t^{\prime})dt^{\prime}%
\label{eq:20241211_142742}
\end{align}
and in position space
\begin{align}
u(x,t)  &  = \iint\frac{ e^{ ik(x-x^{\prime}) } }{2\pi} \frac{e^{-i\omega
_{1}(t+T)}\omega_{2}(k)-e^{-i\omega_{2}(t+T)}\omega_{1}(k)}{\omega_{2}(k)-\omega
_{1}(k)} u(x^{\prime},-T) dk dx^{\prime}\nonumber\\
&  - i \iint\frac{ e^{ ik(x-x^{\prime}) } }{2\pi} \frac{ e^{-i\omega_{1}(t+T)}-e^{-i\omega_{2}(t+T)} }{\omega_{2}(k)-\omega_{1}(k)}\frac{{\partial}
}{{\partial}{t}} u(x^{\prime},t=-T) dk dx^{\prime}\nonumber\\
&  + i\int_{-T}^{0} \left[  \iint\frac{ e^{ ik(x-x^{\prime}) } }{2\pi} \frac{
e^{-i\omega_{1}(t-t^{\prime})} - e^{-i\omega_{2}(t-t^{\prime})} } {\omega
_{2}(k)-\omega_{1}(k)} f(x^{\prime},t^{\prime}) dk dx^{\prime} \right]
dt^{\prime} \label{126_146}%
\end{align}

These equations are the generalization of the d'Alembert solution for the case
of two arbitrary modes; Eq.\ (\ref{126_145}) is in wavenumber space and
Eq.\ (\ref{126_146}) is in position space.

\bigskip

\noindent\textbf{Standard case}.
%
So far,
we have thus obtained the generalization
of the d'Alembert solution
for the case we consider in this paper;
that is,
of a finite source that is on from $-T$ to $0$,
and gave the wave at some later time $t$.
Using our solution,
we obtain the usual case,
which is:
the wave at time $t$
due to a source which is on from $0$ to $t$.
In wavenumber space,
Eq.\ \eqref{126_145}
becomes
\begin{align}
S(k,t)  &  =\frac{e^{-i\omega_{1}t}\omega_{2}(k)-e^{-i\omega_{2}t}\omega
_{1}(k)}{\omega_{2}(k)-\omega_{1}(k)}S(k,0)-i\frac{ e^{-i\omega_{1}%
t}-e^{-i\omega_{2}t} }{\omega_{2}(k)-\omega_{1}(k)}\frac{{\partial}}%
{{\partial}{t}}S(k,0)\nonumber\\
&  +i\int_{0}^{t} \frac{ e^{-i\omega_{1}(t-t^{\prime})} - e^{-i\omega
_{2}(t-t^{\prime})} } {\omega_{2}(k)-\omega_{1}(k)} F(k,t^{\prime})dt^{\prime}
,
\end  {align}
Eq.\ \eqref{eq:20241211_142742}
becomes
\begin{align}
S(k,t)  &  = -i e^{-i(\omega_{2}+\omega_{1})t} \frac{{\partial}}{{\partial}%
{t}} \left(  \frac{ e^{ i\omega_{2}t} - e^{ i\omega_{1}t} } {\omega
_{2}(k)-\omega_{1}(k)} \right)  S(k,0)dk\nonumber\\
&  +i\int_{0}^{t} \frac{ e^{-i\omega_{1}(t-t^{\prime})} - e^{-i\omega
_{2}(t-t^{\prime})} } {\omega_{2}(k)-\omega_{1}(k)} F(k,t^{\prime})dt^{\prime}%
%
\end  {align}
and in position space,
Eq.\ \eqref{126_146}
becomes
\begin{align}
u(x,t)  &  = \iint\frac{ e^{ ik(x-x^{\prime}) } }{2\pi} \frac{e^{-i\omega
_{1}t}\omega_{2}(k)-e^{-i\omega_{2}t}\omega_{1}(k)}{\omega_{2}(k)-\omega
_{1}(k)} u(x^{\prime},0) dk dx^{\prime}\nonumber\\
&  - i \iint\frac{ e^{ ik(x-x^{\prime}) } }{2\pi} \frac{ e^{-i\omega_{1}
t}-e^{-i\omega_{2}t} }{\omega_{2}(k)-\omega_{1}(k)}\frac{{\partial}
}{{\partial}{t}} u(x^{\prime},t=0) dk dx^{\prime}\nonumber\\
&  + i\int_{0}^{t} \left[  \iint\frac{ e^{ ik(x-x^{\prime}) } }{2\pi} \frac{
e^{-i\omega_{1}(t-t^{\prime})} - e^{-i\omega_{2}(t-t^{\prime})} } {\omega
_{2}(k)-\omega_{1}(k)} f(x^{\prime},t^{\prime}) dk dx^{\prime} \right]
dt^{\prime}
%
\end  {align}
%

\section{The $M$-mode case}
\label{sec:20240911_151904}

Suppose that the wave equation produces a dispersion relation that has $M$
modes and that we are given $M$ initial conditions that consist of the wave
$u(x,t)$ and its $M-1$ time-derivatives evaluated at some initial time ($t=-T$).
In terms of modes,
\begin{align}
u(x,-T)  &  =\sum_{\ell=1}^{M}u_{\ell}(x,-T)\\
u{}^{(j)}(x,-T)  &  =\sum_{\ell=1}^{M}\frac{{\partial}^{j}{}}{{\partial}%
{t}^{j}}u_{\ell}(x,-T)\hspace{0.5in}j=1,\cdots\!,M-1
\end{align}
In the wavenumber domain, these are
\begin{align}
S(k,-T)  &  = \sum_{\ell=1}^{M}S_{\ell}(k,-T)\label{eq:20230320_135126}\\
\frac{\partial^{j}}{\partial t^{j}}S(k,-T)  &  = \sum_{\ell=1}^{M}%
\frac{\partial^{j}}{\partial t^{j}}S_{\ell}(k,-T)=\sum_{\ell=1}^{M}\left[
-i\omega_{\ell}\left(  \mathcal{K}\right)  \right]  ^{j}S_{\ell}(k,-T)
\qquad\text{for } j=1,2,\cdots\!,M-1 \label{eq:20230320_135155}%
\end{align}
Eq.\ \eqref{eq:20230320_135126} and \eqref{eq:20230320_135155} are $M$
algebraic equations for the modes $S_{\ell}$. These only hold when the source
is turned off, including at times $t=-T$ and $t=0$. Writing the equations in matrix
form, we have
\begin{equation}%
\begin{pmatrix}
S^{B}(k,-T)\\
\frac{\partial}{\partial t}S^{B}(k,-T)\\
\vdots\\
\frac{\partial^{M-1}}{\partial t^{M-1}}S^{B}(k,-T)
\end{pmatrix}
=%
\begin{pmatrix}
1 & 1 & \cdots & 1\\
-i\omega_{1}(k) & -i\omega_{2}(k) & \cdots & -i\omega_{M}(k)\\
\vdots & \vdots & \ddots & \vdots\\
\lbrack-i\omega_{1}(k)]^{M-1} & [-i\omega_{2}(k)]^{M-1} & \cdots &
[-i\omega_{M}(k)]^{M-1}%
\end{pmatrix}%
\begin{pmatrix}
S_{1}^{B}(k,-T)\\
S_{2}^{B}(k,-T)\\
\vdots\\
S_{M}^{B}(k,-T)
\end{pmatrix}
\end{equation}

We can obtain the individual modes by inverting the matrix,
\begin{equation}%
\begin{pmatrix}
S_{1}^{B}(k,-T)\\
S_{2}^{B}(k,-T)\\
\vdots\\
S_{M}^{B}(k,-T)
\end{pmatrix}
=%
\begin{pmatrix}
1 & 1 & \cdots & 1\\
-i\omega_{1}(k) & -i\omega_{2}(k) & \cdots & -i\omega_{M}(k)\\
\vdots & \vdots & \ddots & \vdots\\
\lbrack-i\omega_{1}(k)]^{M-1} & [-i\omega_{2}(k)]^{M-1} & \cdots &
[-i\omega_{M}(k)]^{M-1}%
\end{pmatrix}
^{-1}%
\begin{pmatrix}
S^{B}(k,-T)\\
\frac{\partial}{\partial t}S^{B}(k,-T)\\
\vdots\\
\frac{\partial^{M-1}}{\partial t^{M-1}}S^{B}(k,-T)
\end{pmatrix}
\end{equation}
Therefore, to obtain the wave some amount of time $\tau$ before time $-T$, we first obtain
the individual modes at time $-T$ from the initial conditions, and propagate
them backward in time by the amount $\tau$; the wave is their superposition
\begin{equation}
S(k,-\tau-T)=\sum_{\ell=1}^{M}e^{i\omega_{\ell}\tau}S_{\ell}^{B}(k,-T),
\end{equation}
which, in matrix form is
\begin{equation}
S(k,-\tau-T)=%
\begin{pmatrix}
e^{i\omega_{1}\tau}\\
e^{i\omega_{2}\tau}\\
\vdots\\
e^{i\omega_{M}\tau}\\
\end{pmatrix}%
\begin{pmatrix}
1 & 1 & \cdots & 1\\
-i\omega_{1}(k) & -i\omega_{2}(k) & \cdots & -i\omega_{M}(k)\\
\vdots & \vdots & \ddots & \vdots\\
\lbrack-i\omega_{1}(k)]^{M-1} & [-i\omega_{2}(k)]^{M-1} & \cdots &
[-i\omega_{M}(k)]^{M-1}%
\end{pmatrix}
^{-1}%
\begin{pmatrix}
S^{B}(k,-T)\\
\frac{\partial}{\partial t}S^{B}(k,-T)\\
\vdots\\
\frac{\partial^{M-1}}{\partial t^{M-1}}S^{B}(k,-T)
\end{pmatrix}
\end{equation}
At time $t=0$, we could obtain the modes from the wave $\varphi$,
\begin{equation}%
\begin{pmatrix}
S_{1}^{A}(k,0)\\
S_{2}^{A}(k,0)\\
\vdots\\
S_{M}^{A}(k,0)
\end{pmatrix}
=%
\begin{pmatrix}
1 & 1 & \cdots & 1\\
-i\omega_{1}(k) & -i\omega_{2}(k) & \cdots & -i\omega_{M}(k)\\
\vdots & \vdots & \ddots & \vdots\\
\lbrack-i\omega_{1}(k)]^{M-1} & [-i\omega_{2}(k)]^{M-1} & \cdots &
[-i\omega_{M}(k)]^{M-1}%
\end{pmatrix}
^{-1}%
\begin{pmatrix}
\varphi(k,0)\\
\frac{\partial}{\partial t}\varphi(k,0)\\
\vdots\\
\frac{\partial^{M-1}}{\partial t^{M-1}}\varphi(k,0)
\end{pmatrix}
\label{eq:20240917_153315}%
\end{equation}
The full wave at some time $t>0$ is
\begin{equation}
S^{A}(k,t)=\sum_{\ell=1}^{M}e^{-i\omega_{\ell}t}S_{\ell}^{A}(k,0)
\end{equation}
which could be expressed in matrix form as
\begin{equation}
S^{A}(k,t)=%
\begin{pmatrix}
e^{-i\omega_{1}t}\\
e^{-i\omega_{2}t}\\
\vdots\\
e^{-i\omega_{M}t}\\
\end{pmatrix}%
\begin{pmatrix}
1 & 1 & \cdots & 1\\
-i\omega_{1}(k) & -i\omega_{2}(k) & \cdots & -i\omega_{M}(k)\\
\vdots & \vdots & \ddots & \vdots\\
\lbrack-i\omega_{1}(k)]^{M-1} & [-i\omega_{2}(k)]^{M-1} & \cdots &
[-i\omega_{M}(k)]^{M-1}%
\end{pmatrix}
^{-1}%
\begin{pmatrix}
\varphi(k,0)\\
\frac{\partial}{\partial t}\varphi(k,0)\\
\vdots\\
\frac{\partial^{M-1}}{\partial t^{M-1}}\varphi(k,0)
\end{pmatrix}
\end{equation}

\bigskip

\subsection{Spatial wave, the $M$-mode case}

Similarly, the spatial modes could be obtained from the initial conditions at
time $-T$
\begin{equation}%
\begin{pmatrix}
u_{1}^{B}(x,-T)\\
u_{2}^{B}(x,-T)\\
\vdots\\
u_{M}^{B}(x,-T)
\end{pmatrix}
=%
\begin{pmatrix}
1 & 1 & \cdots & 1\\
-i\omega_{1}(\mathcal{K}) & -i\omega_{2}(\mathcal{K}) & \cdots & -i\omega
_{M}(\mathcal{K})\\
\vdots & \vdots & \ddots & \vdots\\
\lbrack-i\omega_{1}(\mathcal{K})]^{M-1} & [-i\omega_{2}(\mathcal{K})]^{M-1} &
\cdots & [-i\omega_{M}(\mathcal{K})]^{M-1}%
\end{pmatrix}
^{-1}%
\begin{pmatrix}
u^{B}(k,-T)\\
\frac{\partial}{\partial t} u^{B}(k,-T)\\
\vdots\\
\frac{\partial^{M-1}}{\partial t^{M-1}} u^{B}(k,-T)
\end{pmatrix}
\end{equation}
We can propagate them backward in time to some time $t=-\tau-T$
\begin{equation}
u^{B}(x^{\prime},-\tau-T)=%
\begin{pmatrix}
e^{i\omega_{1}(\mathcal{K})\tau}\\
e^{i\omega_{2}(\mathcal{K})\tau}\\
\vdots\\
e^{i\omega_{M}(\mathcal{K})\tau}\\
\end{pmatrix}%
\begin{pmatrix}
1 & 1 & \cdots & 1\\
-i\omega_{1}(\mathcal{K}) & -i\omega_{2}(\mathcal{K}) & \cdots & -i\omega
_{M}(\mathcal{K})\\
\vdots & \vdots & \ddots & \vdots\\
\lbrack-i\omega_{1}(\mathcal{K})]^{M-1} & [-i\omega_{2}(\mathcal{K})]^{M-1} &
\cdots & [-i\omega_{M}(\mathcal{K})]^{M-1}%
\end{pmatrix}
^{-1}%
\begin{pmatrix}
u^{B}(x,-T)\\
\frac{\partial}{\partial t} u^{B}(x,-T)\\
\vdots\\
\frac{\partial^{M-1}}{\partial t^{M-1}} u^{B}(x,-T)
\end{pmatrix}
\end{equation}

At time $0$, the spatial modes are
\begin{align}%
\begin{pmatrix}
u_{1}(x,0)\\
u_{2}(x,0)\\
\vdots\\
u_{M}(x,0)
\end{pmatrix}
=
\begin{pmatrix}
1 & 1 & \cdots & 1\\
-i\omega_{1}(\mathcal{K}) & -i\omega_{2}(\mathcal{K}) & \cdots & -i\omega
_{M}(\mathcal{K})\\
\vdots & \vdots & \ddots & \vdots\\
\lbrack-i\omega_{1}(\mathcal{K})]^{M-1} & [-i\omega_{2}(\mathcal{K})]^{M-1} &
\cdots & [-i\omega_{M}(\mathcal{K})]^{M-1}%
\end{pmatrix}
^{-1}%
\begin{pmatrix}
\mu(x,0)\\
\frac{\partial}{\partial t}\mu(x,0)\\
\vdots\\
\frac{\partial^{M-1}}{\partial t^{M-1}}\mu(x,0)
\end{pmatrix}
\label{eq:20240917_153426}%
\end{align}
and in position space,
the wave after time $t=0$ is
\begin{equation}
u^{A}(x^{\prime},t)=\iint\frac{dxdk}{2\pi}e^{ik(x^{\prime}-x)}%
\begin{pmatrix}
e^{-i\omega_{1}(\mathcal K)t}\\
e^{-i\omega_{2}(\mathcal K)t}\\
\vdots\\
e^{-i\omega_{M}(\mathcal K)t}\\
\end{pmatrix}
\cdot%
\begin{pmatrix}
u_{1}(x,0)\\
u_{2}(x,0)\\
\vdots\\
u_{M}(x,0)
\end{pmatrix}
\end{equation}
or
\begin{align}
u(x,t) = \sum_{\ell=1}^{M} e^{-i\omega_{\ell}(\mathcal{K})t} u^{A}_{\ell}(x,0)
\end{align}
which may be expressed in matrix form as
\begin{equation}
u^{A}(x,t)=%
\begin{pmatrix}
e^{-i\omega_{1}(\mathcal{K})t}\\
e^{-i\omega_{2}(\mathcal{K})t}\\
\vdots\\
e^{-i\omega_{M}(\mathcal{K})t}\\
\end{pmatrix}%
\begin{pmatrix}
1 & 1 & \cdots & 1\\
-i\omega_{1}(\mathcal{K}) & -i\omega_{2}(\mathcal{K}) & \cdots & -i\omega
_{M}(\mathcal{K})\\
\vdots & \vdots & \ddots & \vdots\\
\lbrack-i\omega_{1}(\mathcal{K})]^{M-1} & [-i\omega_{2}(\mathcal{K})]^{M-1} &
\cdots & [-i\omega_{M}(\mathcal{K})]^{M-1}%
\end{pmatrix}
^{-1}%
\begin{pmatrix}
\mu(x,0)\\
\frac{\partial}{\partial t}\mu(x,0)\\
\vdots\\
\frac{\partial^{M-1}}{\partial t^{M-1}}\mu(x,0)
\end{pmatrix}
\end{equation}

%
%
%

\bigskip
\noindent
\textbf{Acknowledgement.}
%
We gratefully acknowledge the support of the Robert A. Welch Foundation (Grant No. A-1261), Air Force Office of Scientific Research (Award No. FA9550-20-1-0366), and U.S. Department of Energy (Award Number DE-SC-0023103).



\end{document}